\documentclass[english,russian,11pt,twoside,fleqn]{article}
\usepackage{amsmath,amssymb,theorem}
\usepackage[T2A]{fontenc}      
\usepackage[cp1251]{inputenc}  
\usepackage[russianb]{babel} 
\usepackage{multicol}
\usepackage{graphicx}

\theoremstyle{plain}

\begin{document}

\begin{center}
{\bf \Large The Fourier transform and convolutions generated by a
differential operator with boundary condition on a segment}
\end{center}

\vspace{5mm}

\begin{center}
{\bf Baltabek Kanguzhin and Niyaz Tokmagambetov}
\end{center}

\vspace{1mm}

\begin{center}
\emph{Department of Fundamental Mathematics,\\ Faculty of
Mechanics and Mathematics,\\ Al--Farabi Kazakh National
University, ave. Al--Farabi 71, 500040,\\ Almaty, Kazakhstan}
\end{center}

\begin{center}
E--mail: kanbalta@mail.ru, tokmagam@list.ru
\end{center}

\vspace{4mm}

{\bf Abstract.} We introduce the concepts of the Fourier transform
and convolution generated by an arbitrary restriction of the
differentiation operator in the space $L_{2}(0,b).$ In contrast to
the classical convolution, the introduced convolution explicitly
depends on the boundary condition that defines the domain of the
operator $L.$ The convolution is closely connected to the inverse
operator or to the resolvent. So, we first find a representation
for the resolvent, and then introduce the required convolution.

\vspace{4mm}

 {\bf keywords:} fourier transform, convolution, differential
 operator, non--local boundary condition, resolvent, spectrum, coefficient
 functional, basis

\vspace{4mm}

{\bf MSC}:{AMS Mathematics Subject Classification (2000) numbers:
34B10, 34L10, 47G30, 47E05}

\vspace{4mm}

\begin{center}
{\bf 0. Introduction}
\end{center}

The standard Fourier transform is a unitary transform in the
Hilbert space $L_{2}(-\infty, \, +\infty)$ and it is generated by
the operator of differentiation $(-i\frac{d}{dx}),$ because the
system of exponents $\{\exp(i\lambda x), \, \lambda\in R\}$ is a
system of "eigenfunctions" corresponding to its continuous
spectrum. The Fourier transform is closely connected to the
bilinear, commutative, associative convolution without
annihilators. An important fact is that the convolution with the
fundamental solution allows us to find solutions of the
inhomogeneous differential equation, which commutes with
differentiation. Corresponding constructions can be generalized to
arbitrary self-adjoint operators. Instead of the differential
operator $(-i\frac{d}{dx})$ in the space $L_{2}(-\infty, \,
+\infty),$ consider an operator $L$ in the Hilbert space
$L_{2}(0,b),$ where $b<\infty,$ which is generated by the
differential operator $(-i\frac{d}{dx})$ and a boundary condition.
We introduce the concepts of the Fourier transform and convolution
generated by an arbitrary restriction of the differentiation
operator in the space $L_{2}(0,b).$ In contrast to the classical
convolution, the introduced convolution explicitly depends on the
boundary condition that defines the domain of the operator $L.$

As noted above, the convolution is closely connected to the
inverse operator or to the resolvent. So, we first find a
representation for the resolvent, and then introduce the required
convolution.

\vspace{8mm}
\begin{center}
{\bf 1. Resolvent and spectrum of the operator $L$}
\end{center}

Without loss of generality, we assume that the origin belongs to
the resolvent set of the operator $L$, that is, there is an
inverse operator $L^{-1}.$ By M. Otelbaev's theorem [3] such
operators are parameterized by a "boundary" function $\sigma(x)$
from the space $L_{2}(0,b).$

{\bf Theorem 1.1} Let the action of the linear operator $L$ in
$L_{2}(0,b)$ be defined by formula $Ly=-i\frac{dy}{dx}$ with some
(fixed) boundary condition. Suppose there exists the inverse
operator $L^{-1}$ in $L_{2}(0,b).$ Then there is a unique function
$\sigma(x)\in L_{2}(0,b)$ such that the domain of operator $L$ is

$$D(L)=\{y\in W^{1}_{2}[0,b]:  y(0)-\int^{b}_{0}(-i\frac{dy}{dx})\overline{\sigma(x)}dx=0\}.$$

{\bf Proof.} Let us consider equation $Ly=f,$ where $f\in
L_{2}(0,b)$. Since there is the bounded inverse operator $L^{-1}$.
We have $y=L^{-1}f$. Each solution of the differential equation
$-i\frac{dy}{dx}=f$ has the form $y=c+K^{-1}f$, where $c$ is an
arbitrary constant and the operator $K$ corresponds to the Cauchy
problem with zero condition at zero:
$$Ky=-i\frac{dy}{dx},  D(K)=\{y\in W^{1}_{2}[0,b]:  y(0)=0\}.$$

Therefore, the constant $c=L^{-1}f-K^{-1}f$ is dependent on $f$
and represents the value of a bounded linear functional on the
Hilbert space $L_{2}(0,b)$. Then $c=c(f),$ and by the Riesz
theorem on bounded linear functionals on $L_{2}(0,b),$ we have
$$c=\int_{0}^{b}f(x)\overline{\sigma(x)}dx, \sigma(x)\in L_{2}(0,b).$$
The element $\sigma(x)$ is uniquely determined. So, the solutions
of the operator equation $Ly=f$ have the form
$y=\int_{0}^{b}f(x)\overline{\sigma(x)}dx+(K^{-1}f)(x).$ In the
last equation, we substitute $x=0$. As a result, we get
$$y(0)=\int_{0}^{b}f(x)\overline{\sigma(x)}dx.$$

The converse assertion is also true and can be verified directly.

If the function $y(x)$ from $W^{1}_{2}[0,b]$ satisfies condition
$$U(y)=0,$$
then it will belong to the domain $D_{L}$ of $L$, where
$$U(y):=y(0)-\int^{b}_{0}(-i\frac{dy}{dx})\overline{\sigma(x)}dx.$$

Let us denote by $\Delta(\lambda)$ the entire function
$\Delta(\lambda)=1-\lambda\int^{b}_{0}exp(i\lambda
x)\overline{\sigma(x)}dx.$ Then the resolvent of the operator $L$
is
$$
(L-\lambda I)^{-1}f=i\int_{0}^{x}exp(i\lambda(x-\xi))f(\xi)d\xi+
$$
\begin{equation}
\label{Res} +\frac{exp(i\lambda
x)}{\Delta(\lambda)}(\int_{0}^{b}f(x)\overline{\sigma(x)}dx+\lambda
i\int_{0}^{b}\overline{\sigma(x)}dx\int_{0}^{x}exp(i\lambda
(x-\xi))f(\xi)d\xi).
\end{equation}
Indeed, denoting the right--hand side of this equality by $y(x),$
we find it is a direct consequence of
$$y'(x)=if(x)+i\lambda y(x).$$

Let us calculate
$$U(y)=y(0)-\int^{b}_{0}(-i\frac{dy}{dx})\overline{\sigma(x)}dx=
y(0)-\lambda
\int^{b}_{0}y(x)\overline{\sigma(x)}dx-\int^{b}_{0}f(x)\overline{\sigma(x)}dx=$$
$$=\frac{1}{\Delta(\lambda)}(\int_{0}^{b}f(x)\overline{\sigma(x)}dx+\lambda
i\int_{0}^{b}\overline{\sigma(x)}dx\int_{0}^{x}exp(i\lambda(x-\xi))f(\xi)d\xi-$$
$$-\lambda\int_{0}^{b}\overline{\sigma(x)}dx[i\int_{0}^{x}exp(i\lambda(x-\xi))f(\xi)d\xi+$$
$$+\frac{exp(i\lambda x)}{\Delta(\lambda)}(\int_{0}^{b}f(x)\overline{\sigma(x)}dx+\lambda i\int_{0}^{b}\overline{\sigma(x)}dx\int_{0}^{x}exp(i\lambda(x-\xi))f(\xi)d\xi]-\int^{b}_{0}f(x)\overline{\sigma(x)}dx=$$
$$=\frac{1}{\Delta(\lambda)}\int^{b}_{0}f(x)\overline{\sigma(x)}dx+\frac{i\lambda}{\Delta(\lambda)}\int^{b}_{0}\overline{\sigma(x)}dx
\int_{0}^{x}exp(i\lambda(x-\xi))f(\xi)d\xi-$$
$$-\frac{i\lambda}{\Delta(\lambda)}\int^{b}_{0}\overline{\sigma(x)}dx\int_{0}^{x}exp(i\lambda(x-\xi))f(\xi)d\xi+$$
$$+\lambda^{2}i\int^{b}_{0}\overline{\sigma(x)}dx
\int_{0}^{x}exp(i\lambda(x-\xi))f(\xi)d\xi\int^{b}_{0}exp(i\lambda
\mu)\overline{\sigma(\mu)}d\mu-$$
$$-\lambda\frac{exp(i\lambda x)}{\Delta(\lambda)}\int_{0}^{b}\overline{\sigma(x)}dx\int_{0}^{x}f(\xi)\overline{\sigma(\xi)}d\xi-$$
$$-\lambda^{2}i\frac{exp(i\lambda x)}{\Delta(\lambda)}\int_{0}^{b}\overline{\sigma(x)}dx\int_{0}^{b}\overline{\sigma(\mu)}d\mu
\int_{0}^{\mu}exp(i\lambda(\mu-\xi))f(\xi)d\xi-\frac{1}{\Delta(\lambda)}\int^{b}_{0}f(x)\overline{\sigma(x)}dx+$$
$$+\frac{\lambda}{\Delta(\lambda)}\int^{b}_{0}f(\xi)\overline{\sigma(\xi)}d\xi\int^{b}_{0}exp(i\lambda
x)\overline{\sigma(x)}dx=0.$$

The proof is complete.

From the representation of the resolvent by definition of the
spectrum of operator $L$ we get the following theorem.

{\bf Theorem 1.2} The set of zeros with multiplicities of the
entire function $\Delta(\lambda)$ is exactly the same as the
spectrum of the operator $L$.

Since $\Delta(\lambda)$ is the entire function in $\lambda$, we
have that the spectrum of the operator $L$ consists of isolated
eigenvalues of finite multiplicity, and limit points of the
spectrum can only be infinity. From the Paley -- Wiener theorem it
immediately follows that:

{\bf Theorem 1.3} The operator $L$ has either countable number of
eigenvalues, or they are absent. The spectrum is empty if and only
if there exists constant $c\in [0,b]$ such that $\sigma(x)=i$ for
$0\leq x\leq c$ and $\sigma(x)=0$ for $c\leq x\leq b.$

The proof immediately follows from the representation
$$\Delta(\lambda)=exp(i\lambda c)-\int^{c}_{0}exp(i\lambda
x)\overline{(\sigma(x)-i)}dx-\lambda \int^{b}_{c}exp(i\lambda
x)\overline{\sigma(x)}dx$$ and from the Paley--Wiener theorem,
since the presence of the integral term in the right side of this
equation leads to the existence of the growing product
$exp(i\lambda c)\cdot \Delta(\lambda)$.

In what follows, suppose that the indicator diagram of the entire
function $\Delta(\lambda)$ is the segment $[0,ib]$. Then the
spectrum of the operator $L$ is a countable set.

In order for the indicator diagram $\Delta(\lambda)$ to be
represented by the interval $[0, ib],$ it is necessary and
sufficient that
\begin{equation}
\label{ConMinMax} \min(supp(\sigma(x)-i))=0,
\max(supp(\sigma(x)))=b,
\end{equation}
where $supp(g)$ is the support of $g$.

The following theorem is proved just as in the work E.Titmarsh
[4].

{\bf Theorem 1.4} Let the condition ~(\ref{ConMinMax}) hold. Then
the number of zeros $N(r)$ of the function $\Delta(\lambda)$,
which satisfy the inequality $|\lambda|\leq r$, satisfies the
limit inequality
$$lim_{r\rightarrow \infty}\frac{N(r)}{r}=\frac{b}{\pi}.$$

In the work of M. Cartwright [5] it is shown that if $\sigma(x)$
satisfies ~(\ref{ConMinMax}) and it is a function of a bounded
variation, then all the zeros of the function $\Delta(\lambda) $
are in a horizontal strip and the value of $N\frac{br}{\pi}(r)$ is
bounded uniformly by $r$.

\vspace{5mm}
\begin{center}
{\bf 2. Convolution, generated by operator $L$}
\end{center}

To obtain a convolution, we rewrite the resolvent $(L-\lambda
I)^{-1}$ in the form
$$
(L-\lambda I)^{-1}f=i\int_{0}^{x}\frac{exp(i\lambda
(x-\xi))}{\Delta(\lambda)}f(\xi)d\xi
$$
\begin{equation}
\label{Res2}
-\int_{0}^{b}\overline{\sigma(\mu)}d\mu\frac{\partial}{\partial\mu}\Big(
\int_{\mu}^{x}\frac{exp(i\lambda
(x-\xi+\mu))}{\Delta(\lambda)}f(\xi)d\xi\Big).
\end{equation}
     Indeed, from the previous view of the resolvent ~(\ref{Res}),
    we have the chain of equalities:
    $$
    (L-\lambda I)^{-1}f= i\int_{0}^{b}exp(i\lambda
     (x-\xi))f(\xi)d\xi+
     $$

     $$
     +\frac{exp(i\lambda
     x)}{\Delta(\lambda)}(\int_{0}^{b}f(x)\overline{\sigma(x)}dx+\lambda
     i\int_{0}^{b}\overline{\sigma(x)}dx\int_{0}^{x}exp(i\lambda
     (x-\xi))f(\xi)d\xi)=
     $$

$$
=i\int_{0}^{x}\frac{exp(i\lambda
     (x-\xi))}{\Delta(\lambda)}f(\xi)d\xi-
     \int_{0}^{b}\overline{\sigma(\mu)}d\mu\int_{0}^{x}\frac{exp(i\lambda
     (x-\xi+\mu))}{\Delta(\lambda)}f(\xi)d\xi+
     $$

$$
+\int_{0}^{b}\frac{exp(i\lambda
     x)}{\Delta(\lambda)}f(\xi)\overline{\sigma(\xi)}d\xi+\lambda i\int_{0}^{b}\overline{\sigma(\mu)}d\mu\int_{0}^{\mu}\frac{exp(i\lambda
     (x-\xi+\mu))}{\Delta(\lambda)}f(\xi)d\xi=
     $$

     $$
     =i\int_{0}^{x}\frac{exp(i\lambda
     (x-\xi))}{\Delta(\lambda)}f(\xi)d\xi-
     \int_{0}^{b}\overline{\sigma(\mu)}d\mu\int_{\mu}^{x}\frac{exp(i\lambda
     (x-\xi+\mu))}{\Delta(\lambda)}f(\xi)d\xi+
     $$

    $$
    +\int_{0}^{b}\frac{exp(i\lambda
     x)}{\Delta(\lambda)}f(\xi)\overline{\sigma(\xi)}d\xi=
     $$

      $$
      =i\int_{0}^{x}\frac{exp(i\lambda
     (x-\xi))}{\Delta(\lambda)}f(\xi)d\xi-
     \int_{0}^{b}\overline{\sigma(\mu)}d\mu\frac{\partial}{\partial\mu}(
     \int_{\mu}^{x}\frac{exp(i\lambda
     (x-\xi+\mu))}{\Delta(\lambda)}f(\xi)d\xi)
     $$
     This proves the lemma.

{\bf Lemma 2.1} The convolution, defined by the formula
$$
(g\ast
f)(x):=i\int_{0}^{x}g(x-\xi)f(\xi)d\xi-\int_{0}^{b}\overline{\sigma(\mu)}d\mu\frac{\partial}{\partial\mu}(
     \int_{\mu}^{x}g(x-\xi+\mu)f(\xi)d\xi)
$$
for $g,f\in W^{1}_{2}[0,b]$ is bilinear, commutative and
associative.

{\bf Proof.} Let us introduce an operation $\circ$ as
$$
(g\circ f)(x,t):=\int_{t}^{x}g(x-\xi+\mu)f(\xi)d\xi,
$$
then we can rewrite the expression $g\ast f$ in the form
$$
(g\ast f)=i(g\circ
f)(x,0)-\int_{0}^{b}\frac{\partial}{\partial\mu}(g\circ
f)(x,\mu)\overline{\sigma(\mu)}d\mu=
$$
$$
=i[(g\circ
f)(x,0)-\int_{0}^{b}(-i\frac{\partial}{\partial\mu}(g\circ
f)(x,\mu))\overline{\sigma_{0}(\mu)}d\mu]=iU_{t}(g\circ f)(x,t),
$$
where
$$
U(y)=y(0)-\int_{0}^{b}(-i\frac{\partial}{\partial\mu})\overline{\sigma(\mu)}d\mu.
$$

Hence, bilinearity is obvious. It is easy to see that the
operation $\circ$ is commutative. Indeed
$$
(g\circ f)(x,t)=\int_{t}^{x}g(x+t-\xi)f(\xi)d\xi=
\left(%
\begin{array}{c}
  \xi'=x+t-\xi \\
  d\xi=-d\xi' \\
  x\rightarrow t\\
 t\rightarrow x \\
\end{array}%
\right)=$$
$$=-\int_{x}^{t}g(\xi')f(x+t-\xi')d\xi'=\int_{t}^{x}f(x+t-\xi)g(\xi)d\xi=(f\circ
g)(x,t)$$
i.е.
$$
(g\circ f)(x,t)=(f\circ g)(x,t)
$$
Then the operation $\ast$ commutative too,
$$
(g\ast f)(x)=i(g\circ
f)(x,0)-\int_{0}^{b}\frac{\partial}{\partial\mu}(g\circ
f)(x,\mu)\overline{\sigma(\mu)}d\mu=
$$
$$
=i(f\circ g)(x,0)-\int_{0}^{b}\frac{\partial}{\partial\mu}(f\circ
g)(x,\mu)\overline{\sigma(\mu)}d\mu=(f\ast g)(x)
$$

{\bf Remark 2.1} The convolution is expressed in terms of boundary
conditions and has the form
\begin{equation}
\label{Conv3} (f\ast
g)(x)=iU_{\mu}\{\int_{\mu}^{x}f(\xi)g(x+\mu-\xi)d\xi\},
\end{equation}
where
$$
U(y):=y(0)-\int^{b}_{0}(-i\frac{dy}{dx})\overline{\sigma(x)}dx.
$$

{\bf Remark 2.2} From ~(\ref{Res2}) it is easy to see that the
resolvent in terms of convolution is represented in the form
\begin{equation}
\label{Res3} (L-\lambda I)^{-1}f=\frac{exp(i\lambda
x)}{\Delta(\lambda)}\ast f(x).
\end{equation}

{\bf Lemma 2.2} For any $f$  from the domain $D_{L}$ and for
arbitrary $g$ from $W^{1}_{2}[0,b],$ the equality
$$
\frac{d}{dx}(f\ast g)=\frac{df}{dx}\ast g
$$
holds true.

{\bf Proof.} Since the equality $U(f)=0$ holds for $f\in D_{L}$,
we get
$$\frac{d}{dx}(f\ast
 g)=if(0)g(x)+i\int_{0}^{x}\frac{df(x-\xi)}{dx}g(\xi)d\xi-$$
 $$
 -\int^{b}_{0}\overline{\sigma(\mu)}d\mu\frac{\partial}{\partial\mu}\begin{pmatrix}
   f(\mu)g(x)+\int_{\mu}^{x}\frac{\partial}{\partial x}f(x-\xi+\mu)g(\xi)d\xi \\
 \end{pmatrix}=\frac{df}{dx}\ast g+if(0)g(x)-
 $$
 $$
 -\int^{b}_{0}\overline{\sigma_{0}(x)}\frac{\partial}{\partial\mu}f(\mu)g(x)d\mu=\frac{df}{dx}\ast g+ig(x)U(f)=\frac{df}{dx}\ast g.
 $$

{\bf Lemma 2.3} The convolution defined in Lemma 2.1 is without
annihilators, i.e. if for arbitrary $$g\in W^{1}_{2}[0,b]$$ and
for some $f\in L_{2}(0,b)$ the equality
$$(g\ast f)(x)=0$$
holds, then $f$ is identically equal to zero.

{\bf Proof.} Let $(g\ast f)(x)=0$,   $f\in L_{2}(0,b),$ then we
take
$$
g=:\frac{exp(i\lambda x)}{\Delta(\lambda)}|_{\lambda=0}=1.
$$
By the definition of our convolution, the expression $1\ast f=0$
denotes the equality $L^{-1}f=0$. Denoting $L^{-1}f$ by $y$, we
obtain respectively $Ly =f,$ but $y =0$, hence, as the $L$ a
linear operator, it follows that $f=0$.

{\bf Lemma 2.4} If the function $f$ is from the domain $D_{L}$ of
the operator $L$, then for each $g\in L_{2}(0,b)$ the convolution
$f\ast g$ will also belong to the domain $D_{L}$.

{\bf Proof.} Let
$$
f=:\frac{exp(i\lambda _{0}x)}{\Delta(\lambda _{0})}|_{\lambda
_{0}=0}=1.
$$
Let us denote by $y$ the function $y=1\ast g$, by construction of
the convolution this expression means that $y=L^{-1}g$ i.e.
$Ly=g$, which implies $y\in D_{L}$. Now we fix any $\lambda$ such
that $\Delta(\lambda)\neq0$. Let $f=:\frac{exp(i\lambda
x)}{\Delta(\lambda)}$.

If by $y$ we denote $y=1\ast g$, then this expression by the
convolution construction mean that $y=L^{-1}g,$ i.e. $Ly=g$, and
implies $y\in D_{L}.$ Now we fix any $\lambda$, such that
$\Delta(\lambda)\neq0$. Let $f=:\frac{exp(i\lambda
x)}{\Delta(\lambda)}$. Let $y$ denote the convolution
$y=\frac{exp(i\lambda x)}{\Delta(\lambda)}\ast g$, then from the
definition of the convolution, we get $y=(L-\lambda I)^{-1}g,$
i.e. $Ly=\lambda y+g$ for any $\lambda$. Since zeros of the
function $\Delta(\lambda)$ are countable set, then there is
sequence $\{\lambda_{n}\}$ such that $\Delta(\lambda_{n})\neq0$
and the system $exp(i\lambda _{n}x)$ is basis in the space
$W^{1}_{2}[0,b]$. For any natural number $n$ it is $\lambda_{n}$
such that the convolution $exp(i\lambda _{n}x)\ast g$ will belong
to the domain $D_{L}$. From
$f(x)=\sum_{n=0}^{\infty}c_{n}exp(i\lambda _{n}x)$ and from the
bilinearity of the convolution, the convolution $f\ast g$ belongs
to the domain $D_{L}$.

{\bf Lemma 2.5} For arbitrary $\lambda$ and $\beta$ such that
$\lambda\neq\beta$ we have the equality
$$
exp(i\lambda x)\ast exp(i\beta x)=\frac{exp(i\beta
x)\Delta(\lambda)-exp(i\lambda x)\Delta(\beta)}{\beta-\lambda}.
$$
{\bf Proof.} We write by definition
$$
exp(i\lambda x)\ast exp(i\beta x)=
$$
$$
=i\int_{0}^{x}exp(i\lambda (x-\xi)) exp(i\beta \xi)d\xi-
$$
$$
-\int^{b}_{0}\overline{\sigma(\mu)}d\mu\frac{\partial}{\partial\mu}
\begin{pmatrix}
\int_{\mu}^{x}exp(i\lambda (x-\xi+\mu)) exp(i\beta \xi)d\xi \\
\end{pmatrix}=
$$
$$
=\frac{iexp(i\lambda
x)(exp(i(\beta-\lambda)x)-1)}{i(\beta-\lambda)}-
$$
$$
-\int^{b}_{0}\overline{\sigma(\mu)}d\mu\frac{\partial}{\partial\mu}
\begin{pmatrix}
exp(i\lambda (x+\mu)) \frac{exp(i(\beta -\lambda)x)-exp(i(\beta -\lambda)\mu)}{i(\beta-\lambda)} \\
\end{pmatrix}=
$$
$$
=\frac{exp(i\beta x)- exp(i\lambda x)}{\beta-\lambda}-
$$
$$
-\int^{b}_{0}\overline{\sigma(\mu)}d\mu\frac{\partial}{\partial\mu}\begin{bmatrix}
\frac{exp(i\lambda \mu )exp(i\beta x)-exp(i\lambda x )exp(i\beta \mu)}{i(\beta-\lambda)} \\
\end{bmatrix}=
$$
$$
=\frac{exp(i\beta x)- exp(i\lambda x)}{\beta-\lambda}-
$$
$$
-\int^{b}_{0}\frac{i\lambda exp(i\lambda \mu )exp(i\beta x)-i\beta
exp(i\lambda x )exp(i\beta
\mu)}{i(\beta-\lambda)}\overline{\sigma(\mu)}d\mu=
$$
$$
\frac{1}{\beta-\lambda}[exp(i\beta x)-exp(i\lambda x )-\lambda
exp(i\beta x)\int^{b}_{0}exp(i\lambda \mu
)\overline{\sigma(\mu)}d\mu+
$$
$$
+\beta exp(i\lambda x )\int^{b}_{0}exp(i\beta
\mu)\overline{\sigma(\mu)}d\mu]=\frac{exp(i\beta
x)\Delta(\lambda)-exp(i\lambda x)\Delta(\beta)}{\beta-\lambda}.
$$

\vspace{5mm}
\begin{center}
{\bf 3. Fourier transform and convolution}
\end{center}

For any function $f$ from the space $L_{2}(0,b),$ let us associate
the expansion
$$
f\sim \sum_{\lambda_{n}\in \sigma(L)}P_{n}f,
$$
where the orthogonal projection is given by
$$
P_{n}f=-\frac{1}{2\pi
i}\oint_{|\lambda-\lambda_{n}|=\delta}(L-\lambda I)^{-1}fd\lambda.
$$
By using ~(\ref{Res3}) we obtain
$$
P_{n}f=res_{\lambda_{n}}\frac{exp(i\lambda
x)}{\Delta(\lambda)}\ast f=
[\frac{1}{(m_{n}-1)!}lim_{\lambda\rightarrow\lambda_{n}}\frac{\partial^{m_{n}-1}}{\partial
\lambda^{m_{n}-1}}(\frac{(\lambda-\lambda_{n})^{m_{n}}}{\Delta(\lambda)}exp(i\lambda
x))]\ast f(x)=
$$
$$
=[\sum_{j=0}^{m_{n}-1}\frac{d_{j,n}}{j!}(\frac{(ix)^{m_{n}-1-j}}{(m_{n}-1-j)!}exp(i\lambda_{n}
x))]\ast f(x).
$$

Let us define
\begin{equation}
\label{RSp}
u_{m_{n}-1,n}:=\sum_{j=0}^{m_{n}-1}\frac{d_{j,n}}{j!}(\frac{(ix)^{m_{n}-1-j}}{(m_{n}-1-j)!}exp(i\lambda_{n}
x)),
\end{equation}
where $\{d_{j,n}\}$ are the expansion coefficients in the Taylor
series of the function
$\frac{(\lambda-\lambda_{n})^{m_{n}}}{\Delta(\lambda)}$ in front
of powers of $(\lambda-\lambda_{n}).$ Then
\begin{equation}
\label{RSp2} P_{n}f=u_{m_{n}-1,n}\ast f.
\end{equation}

Let us introduce a system of root functions of the operator $L$
corresponding to the eigenvalue $\lambda_{n}$ by the expression
\begin{equation}
\label{RSp3}
u_{s,n}:=Lu_{s+1,n}-\lambda_{n}u_{s+1,n}=\sum_{j=0}^{s}\frac{d_{j,n}}{j!}\frac{(ix)^{s-j}}{(s-j)!}exp(i\lambda_{n}
x), \,\,\, s=0,.....,m_{n}-2.
\end{equation}

{\bf Lemma 3.1} The system of root functions $\{u_{k,n},
k=0...,m_{n}-1\}$ is linearly independent.

{\bf Proof.} Consider the linear combination
\begin{equation}
\label{LinComb}
\alpha_{0}u_{0,n}+...+\alpha_{m_{n}-2}u_{m_{n}-2,n}+\alpha_{m_{n}-1}u_{m_{n}-1,n}=0.
\end{equation}

It is easy to see that
\begin{equation}
\label{LinComb2} (L-\lambda_{n} I)^{m_{n}-1}u_{s,n}=0.
s=0,...,m_{n}-2.
\end{equation}
By applying to ~(\ref{LinComb}) the operator
$(L-\lambda_{n}I)^{m_{n}-1}$ on both sides, we get
\begin{equation}
\label{LinComb3}\alpha_{m_{n}-1}(L-\lambda_{n}
I)^{m_{n}-1}u_{m_{n}-1,n}=0.
\end{equation}
From ~(\ref{RSp3}) we see that ~(\ref{LinComb3}) is equivalent to
$$
\alpha_{m_{n}-1}u_{0,n}=0.
$$
Since $d_{0,n}\neq0$, we obtain
$$
\alpha_{m_{n}-1}=0.
$$
Now we consider the equality ~(\ref{LinComb}) with
$\alpha_{m_{n}-1}=0$, i.e.
$$
\alpha_{0}u_{0,n}+...+\alpha_{m_{n}-2}u_{m_{n}-2,n}=0.
$$
Applying to both sides the operator $(L-\lambda_{n}I)^{m_{n}-2}$
and using ~(\ref{LinComb2}) as in the previous step, we get
$$
\alpha_{m_{n}-2}u_{0,n}=0,
$$
whence $\alpha_{m_{n}-2}=0.$

Repeating these steps, we find that equality ~(\ref{LinComb}) is
true only for $\alpha_{k}=0, k=0,...,m_{n}-1$. This proves the
linear independence of the system $\{u_{k,n}, k=0...,m_{n}-1\}$.

Since $d_{0,n}\neq0$, the system
$u_{0,n},u_{1,n},...,u_{m_{n}-1,n}$ is linearly independent and,
therefore, there is a basis of the root subspace
$$
H_{\lambda_{n}}:=Ker(L-\lambda_{n} I)^{m_{n}}.
$$
We expand the function $P_{n}f$ on $Ker(L-\lambda_{n}I)^{m_{n}}$
by this basis
$$
P_{n}f=C_{m_{n}-1,n}(f)u_{0,n}+C_{m_{n}-2,n}(f)u_{1,n}+...+C_{0,n}(f)u_{m_{n}-1,n}.
$$
Thus, for each element $f$ of the space $L_{2}(0, b)$ we associate
the element of the sequence space
$$
\{(C_{0,n}(f),...,C_{m_{n}-1,n}(f)),\lambda_{n}\in \sigma(L)\}\in
\prod_{\lambda_{n}\in \sigma(L)}\mathbb C^{m_{n}},
$$
i.e. we introduce the Fourier transform
$$
\widehat{f}=\{(C_{0,n}(f),...,C_{m_{n}-1,n}(f)),\lambda_{n}\in
\sigma(L)\}.
$$

In the space of sequences, we introduce the inner Cauchy
convolution. Let $\xi$ and $\eta$ be any elements of the
$X:=\prod_{\lambda_{n}\in \sigma(L)}\mathbb C^{m_{n}}$, then we
will call their convolution the sequence
$$
\mu:=\{(\xi_{0,n}\eta_{0,n};\xi_{0,n}\eta_{1,n}+\xi_{1,n}\eta_{0,n};...,\sum_{k=0}^{m_{n}-1}\xi_{k,n}\eta_{m_{n}-1-k,n}),
\lambda_{n}\in \sigma(L)\},
$$
which we denote by $\xi\ast_{X}\eta$. Introduced convolutions
$\ast_{X}$ and $\ast$ are associated between themselves by the
Fourier transform.

{\bf Theorem 3.1} For arbitrary functions $f$ and $g$ from the
space $W^{1}_{2}[0,b]$ the equality
\begin{equation}
\label{FTr} \widehat{f\ast g}=\widehat{f}\ast_{X}\widehat{g}
\end{equation}
holds.

{\bf Proof.} Let $f=u_{m_{n}-1,n}$. Then the equality
$\widehat{u_{m_{n}-1,n}\ast g}=
\widehat{u}_{m_{n}-1,n}\ast_{X}\widehat{g}$ follows from the
following chain of equalities
$$
(u_{m_{n}-1,n}\ast
g)\hat{}=(P_{n}g)\hat{}=(C_{0,n}(g),C_{1,n}(g),...,C_{m_{n}-1,n}(g))=
$$
$$
=(1,0,...,0)\ast_{X}(C_{0,n}(g),C_{1,n}(g),...,C_{m_{n}-1,n}(g)).
$$

Let now $f=u_{s,n}$, then
$$
(u_{s,n}\ast g)\hat{}=((L-\lambda_{n}
I)^{m_{n}-1-s}u_{m_{n}-1,n}\ast g)\hat{}=
$$
$$
=((L-\lambda_{n} I)^{m_{n}-1-s}(u_{m_{n}-1,n}\ast
g))\hat{}=((L-\lambda_{n}
I)^{m_{n}-1-s}P_{n}g))\hat{}=(P_{n}(L-\lambda_{n}
I)^{m_{n}-1-s}g))\hat{}=
$$
$$
=(C_{0,n})(L-\lambda_{n} I)^{m_{n}-1-s}g)),(C_{1,n})(L-\lambda_{n}
I)^{m_{n}-1-s}g)),...,(C_{m_{n}-1,n})(L-\lambda_{n}
I)^{m_{n}-1-s}g)).
$$
On the other hand, the relations
$$
C_{m_{n}-1,n}((L-\lambda_{n} I)g)=C_{m_{n}-2,n}(g),
$$
$$...,$$
$$
C_{1,n}((L-\lambda_{n} I)g)=C_{0,n}(g)
$$
and
$$
C_{0,n}((L-\lambda_{n} I)g)=0
$$
are valid, i.e. the action of the operator $(L-\lambda_{n} I)$ is
equivalent to a shift in the sequence space. Then action of the
operator $(L-\lambda_{n} I)^{m_{n}-1-s}$ corresponds to a shift to
the right on $m_{n}-1-s$ position. This implies the equality
$$
\widehat{u_{s,n}\ast g}=\widehat{u_{s,n}}\ast_{X}\widehat{g},
\,\,\, s=0,...,m_{n}-1.
$$
From the fact that the theorem holds for all elements of the basis
$\{u_{s,n}\}$ follows the required equality ~(\ref{FTr}) follow in
the whole space.

\vspace{5mm}
\begin{center}
{\bf 4. Coefficient functionals and a boundary condition}
\end{center}

Elements of the basis of the each root subspace have good
properties.

{\bf Lemma 4.1} Let $\lambda_{n}$ be zeros of the entire function
$\Delta(\lambda)$ with corresponding multiplicities $m_{n}$. Then
the elements ~(\ref{RSp}) and ~(\ref{RSp3}) of the root subspace
$H_{\lambda_{n}}$ have the following properties:
\begin{equation}
\label{ConvProd} u_{p,n}\ast u_{q,n}=\begin{cases}
0,&\text{for $p+q<m_{n}-1$},\\
u_{p+q-m_{n}+1},&\text{for $p+q\geq m_{n}-1$,}\\
\end{cases}
\end{equation}
$0\leq p,q\leq m_{n}-1.$

{\bf Proof.} At first, we note that the element ~(\ref{RSp}) is
idempotent with respect to the convolution, i.e.
$$
u_{m_{n}-1,n}\ast u_{m_{n}-1,n}=u_{m_{n}-1,n}.
$$

By acting on the element $u_{m_{n}-1, n}$ of the basis of root
subspace of the operator $P_{n}$, we get
$$
P_{n}(u_{m_{n}-1,n})=u_{m_{n}-1,n}.
$$
On the other hand, by replacing the function $f$ by
$u_{m_{n}-1,n}$ in the formula ~(\ref{RSp2}), we have
$$
P_{n}(u_{m_{n}-1,n})=u_{m_{n}-1,n}\ast u_{m_{n}-1,n},
$$
which proves the property that the element $u_{m_{n}-1,n}$ is
idempotent. Further, from ~(\ref{RSp3}) it is easy to see that
each basis element can be presented by an idempotent element
$u_{m_{n}-1,n}$ in the form
\begin{equation}
\label{ConvProd3} u_{q,n}=(L-\lambda_{n}
I)^{m_{n}-1-q}u_{m_{n}-1,n},  q=0,...,m_{n}-2.
\end{equation}
Then consider the convolution of two elements of the basis
$$
u_{q,n}\ast u_{p,n}=[(L-\lambda_{n}
I)^{m_{n}-1-q}u_{m_{n}-1,n}]\ast [(L-\lambda_{n}
I)^{m_{n}-1-p}u_{m_{n}-1,n}]=
$$
|by using bilinear property of the convolution $\ast$|
$$
=(L-\lambda_{n} I)^{2m_{n}-2-q-p}(u_{m_{n}-1,n}\ast
u_{m_{n}-1,n})=
$$
|by using idempotence property of the element $u_{m_{n}-1,n}$|
$$
=(L-\lambda_{n} I)^{m_{n}-1-(q+p-m_{n}+1}u_{m_{n}-1,n}.
$$
It is clear that if the inequality $q+p-m_{n}+1\geq 0$ holds, then
from ~(\ref{ConvProd3}) it follows that
$$u_{p,n}\ast u_{q,n}=u_{p+q-m_{n}+1,n}$$
and otherwise
$$u_{p,n}\ast u_{q,n}=0.$$

In the following lemma we give the coefficient functionals in the
boundary condition.

{\bf Lemma 4.2} The coefficient functionals in the projector
expansion $P_{n}f=\sum_{k=0}^{m_{n}-1}C_{m_{n}-1-k,n}(f)u_{k,n}$
in the root subspace $Ker(L-\lambda_{n} I)^{m_{n}}$ are of the
form
$$
C_{k,n}(f)=-iU_{\mu}\{\int_{0}^{\mu}f(\xi)\frac{(i(\mu-\xi))^{k}}{k!}exp(i\lambda_{n}(\mu-\xi))d\xi\},
0\leq k\leq m_{n}-1.
$$

{\bf Proof.} By the definition we have $P_{n}f=f\ast
u_{m_{n}-1,n}$, but on other hand
$P_{n}f=C_{m_{n}-1,n}(f)u_{0,n}+C_{m_{n}-2,n}(f)u_{1,n}+...+C_{0,n}(f)u_{m_{n}-1,n}.$
Then we get the equality
\begin{equation}
\label{Coef2} f\ast
u_{m_{n}-1,n}=\sum_{k=0}^{m_{n}-1}C_{m_{n}-1-k,n}(f)u_{k,n},
\end{equation}
As for $u_{0,n}$, we have
$$
f\ast u_{m_{n}-1,n}\ast
u_{0,n}=\sum_{k=0}^{m_{n}-1}C_{m_{n}-1-k,n}(f)u_{k,n}\ast u_{0,n},
$$
then by using Lemma 3.1, we get
\begin{equation}
\label{Coef4} f\ast u_{0,n}=C_{0,n}(f)u_{0,n}.
\end{equation}
Let rewrite the left side of ~(\ref{Coef4}) by using formula
~(\ref{Conv3})
$$
f\ast
u_{0,n}=-iU_{\mu}\{\int_{\mu}^{x}f(\xi)u_{0,n}(x+\mu-\xi)d\xi\}=
$$
$$
=-iU_{\mu}\{\int_{0}^{x}f(\xi)u_{0,n}(x+\mu-\xi)d\xi\}-iU_{\mu}\{\int_{0}^{\mu}f(\xi)u_{0,n}(x+\mu-\xi)d\xi\}.
$$
Since
$iU_{\mu}\{\int_{0}^{x}f(\xi)d_{0,n}exp(i\lambda_{n}(x+\mu-\xi))d\xi\}=$
$$
=d_{0,n}exp(i\lambda
x)i\int_{0}^{x}f(\xi)exp(-i\lambda_{n}\xi)d\xi[1-\lambda_{n}\int_{0}^{b}exp(i\lambda_{n}\mu)\overline{\sigma(\mu)}d\mu]=
$$
$$
=u_{0,n}i\Delta(\lambda_{n})\int_{0}^{x}f(\xi)exp(-i\lambda_{n}\xi)d\xi,
$$
the first term is equal to zero. And
$$
f\ast
u_{0,n}=-iU_{\mu}\{\int_{0}^{\mu}f(\xi)u_{0,n}(x+\mu-\xi)d\xi\}=-iU_{\mu}\{\int_{0}^{\mu}f(\xi)d_{0,n}exp(i\lambda_{n}(x+\mu-\xi))d\xi\}=
$$
$$
=-d_{n,0}exp(-i\lambda_{n}\xi)iU_{\mu}\{\int_{0}^{\mu}f(\xi)exp(i\lambda_{n}(\mu-\xi))d\xi\}=-iU_{\mu}\{\int_{0}^{\mu}f(\xi)exp(i\lambda_{n}(\mu-\xi))d\xi\}u_{0,n}.
$$
By comparing the right and the left hand side of ~(\ref{Coef4}),
we get
$$
C_{0,n}(f)=-iU_{\mu}\{\int_{0}^{\mu}f(\xi)exp(i\lambda_{n}(\mu-\xi))d\xi\}.
$$
Now suppose that
$$
C_{k,n}(f)=-iU_{\mu}\{\int_{0}^{\mu}f(\xi)\frac{(i(\mu-\xi))^{k}}{k!}exp(i\lambda_{n}(\mu-\xi))d\xi\},
k=0:s-1,  1\leq s\leq m_{n}-1.
$$
Let us prove that
$$
C_{s,n}(f)=-iU_{\mu}\{\int_{0}^{\mu}f(\xi)\frac{(i(\mu-\xi))^{s}}{s!}exp(i\lambda_{n}(\mu-\xi))d\xi\},
1\leq s\leq m_{n}-1.
$$
By convolving the both sides of ~(\ref{Coef2}) with $u_{s,n}$, we
have
\begin{equation}
\label{Coef8} f\ast u_{m_{n}-1,n}\ast
u_{s,n}=\sum_{k=0}^{m_{n}-1}C_{m_{n}-1-k,n}(f)u_{k,n}\ast u_{s,n}.
\end{equation}
The equality ~(\ref{Coef8}) is changed to
$$
f\ast
u_{s,n}=\sum_{k=m_{n}-1-s}^{m_{n}-1}C_{m_{n}-1-k,n}(f)u_{k+s-m_{n}+1,n}.
$$
Let us related the index $k$ by $l=m_{n}-1-k,$ then we get
\begin{equation}
\label{Coef10} f\ast u_{s,n}=\sum_{l=0}^{s}C_{l,n}(f)u_{s-1,n}.
\end{equation}
As in the previous case by using formula ~(\ref{Conv3}), we write
$$
f\ast
u_{s,n}=iU_{\mu}\{\int_{\mu}^{x}f(\xi)u_{s,n}(x+\mu-\xi)d\xi\}=
$$
$$
=iU_{\mu}\{\int_{0}^{x}f(\xi)u_{s,n}(x+\mu-\xi)d\xi\}-iU_{\mu}\{\int_{0}^{\mu}f(\xi)u_{s,n}(x+\mu-\xi)d\xi\}.
$$
Since
$$
iU_{\mu}\{\int_{0}^{x}f(\xi)\sum_{j=0}^{s}\frac{d_{j,n}}{j!}\frac{(i(x+\mu-\xi))^{s-j}}{(s-j)!}exp(i\lambda_{n}(x+\mu-\xi))d\xi\}=
$$
$$
=\sum_{j=0}^{s}\frac{d_{j,n}}{j!}exp(i\lambda_{n}x)\sum_{p=0}^{s-j}\frac{(ix)^{s-j-p}}{(s-j-p)!}iU_{\mu}
\{\int_{0}^{x}f(\xi)\frac{(i(\mu-\xi))^{p}}{p!}exp(i\lambda_{n}(\mu-\xi))d\xi\}=
$$
$$
=\sum_{j=0}^{s}\frac{\partial^{j}}{\partial\lambda^{j}}G(f;x,\lambda)|_{\lambda=\lambda_{n}}u_{s-j,n}
$$
the first term is equal to zero, here
$$
G(f;x,\lambda):=i\Delta(\lambda)\int_{0}^{x}f(\xi)exp(-i\lambda\xi)d\xi.
$$

It is obvious that if $\lambda_{n}$ is zero of the function
$\Delta(\lambda)$, then
$\frac{\partial^{j}}{\partial\lambda^{j}}G(f;x,\lambda)|_{\lambda=\lambda_{n}}=0,
\forall j\in Z_{+}$.

In the left hand side of ~(\ref{Coef10}), we then have
$$
f\ast
u_{s,n}=-iU_{\mu}\{\int_{0}^{\mu}f(\xi)u_{s,n}(x+\mu-\xi)d\xi\}=
$$
$$
=-iU_{\mu}\{\int_{0}^{\mu}f(\xi)\sum_{j=0}^{s}\frac{d_{j,n}}{j!}\frac{(i(x+\mu-\xi))^{s-j}}{(s-j)!}exp(i\lambda_{n}(x+\mu-\xi))d\xi\}=
$$
$$
=-\sum_{j=0}^{s}\frac{d_{j,n}}{j!}exp(i\lambda_{n}x)\sum_{p=0}^{s-j}\frac{(ix)^{s-j-p}}{(s-j-p)!}iU_{\mu}
\{\int_{0}^{\mu}f(\xi)\frac{(i(\mu-\xi))^{p}}{p!}exp(i\lambda_{n}(\mu-\xi))d\xi\}=
$$
$$
=\sum_{j=0}^{s}C_{j,n}(f)u_{s-j,n}-iU_{\mu}\{\int_{0}^{\mu}f(\xi)\frac{(i(\mu-\xi))^{s}}{s!}exp(i\lambda_{n}(\mu-\xi))d\xi\}u_{0,n}.
$$
By comparing the right and the left hand sides of the equality
~(\ref{Coef10}), we get
$$
C_{s,n}(f)=-iU_{\mu}\{\int_{0}^{\mu}f(\xi)\frac{(i(\mu-\xi))^{s}}{s!}exp(i\lambda_{n}(\mu-\xi))d\xi\}.
$$
Hence we have explicitly constructed a biorthogonal system to the
system $\{u_{k,n},\lambda_{n}\in \sigma(L)\}.$

{\bf Theorem 4.2} The chosen basis system $\{u_{k,n},
k=0,...,m_{n}-1,\lambda_{n}\in \sigma(L)\}$ is minimal in
$L_{2}(0,b),$ i.e. there exists a biorthogonal system of the form
$$
h_{k,n}(\xi)=\int_{\xi}^{b}\overline{\sigma(\mu)}\frac{\partial}{\partial\mu}(\frac{(i(\mu-\xi))^{k}}{k!}exp(i\lambda_{n}(\mu-\xi)))d\mu,
k=0,...,m_{n}-1.
$$

\vspace{5mm}
\begin{center}
{\bf 5. On the Sedletskiy formula for the remainder term}
\end{center}

It is well--known that the partial sum of the Fourier series is
written in the form
$$
S_{R}(f;x)=-\frac{1}{2\pi i}\oint_{|\lambda|=R}(L-\lambda
I)^{-1}fd\lambda=\sum_{|\lambda_{n}|<R} P_{\lambda_{n}}f.
$$
By using ~(\ref{Res3}), we rewrite the partial sum as
$$
S_{R}(f;x)=-\frac{1}{2\pi i}\oint_{|\lambda|=R}\frac{exp(i\lambda
x)}{\Delta (\lambda)}\ast f(x)d\lambda.
$$
In this formula instead of $f$ we put the function $exp(i\mu x).$
Then we have
$$
S_{R}(exp(i\mu x);x)=-\frac{1}{2\pi
i}\oint_{|\lambda|=R}\frac{exp(i\lambda x)}{\Delta (\lambda)}\ast
exp(i\mu x)d\lambda=
$$ |by applying lemma 2.5| $=$
$$
=\frac{1}{2\pi i}\oint_{|\lambda|=R}\frac{exp(i\lambda
x)\Delta(\mu)-exp(i\mu x)\Delta(\lambda)}{\Delta
(\lambda)(\mu-\lambda)}d\lambda=
$$
$$
=\frac{\Delta(\mu)}{2\pi i}\oint_{|\lambda|=R}\frac{exp(i\lambda
x)}{\Delta (\lambda)}\frac{d\lambda}{\mu-\lambda}+\frac{1}{2\pi
i}\oint_{|\lambda|=R}\frac{exp(i\mu x)}{\lambda-\mu}d\lambda=
$$
|by using Cauchy formula to the second term| $=$
$$
=exp(i\mu x)+\frac{\Delta(\mu)}{2\pi
i}\oint_{|\lambda|=R}\frac{exp(i\lambda x)}{\Delta
(\lambda)}\frac{d\lambda}{\mu-\lambda}.
$$
Hence it follows that the formula for the remainder term for
$f(x)=exp(i\mu x)$ is
\begin{equation}
\label{Q} Q_{R}(exp(i\mu x);x):=S_{R}(exp(i\mu x);x)-exp(i\mu
x)=\frac{\Delta(\mu)}{2\pi i}\oint_{|\lambda|=R}\frac{exp(i\lambda
x)}{\Delta (\lambda)}\frac{d\lambda}{\mu-\lambda},
\end{equation}
which is valid for any $\mu$. Since an arbitrary element $f$ of
the space $L_{2}(0,b)$  is represented in the form
$$
f(x)=\frac{1}{2\pi}\int_{-\infty}^{\infty}\hat{f}(-\mu)exp(i\mu
x)d\mu,
$$
where $\widehat{f}(\mu)$ is the Fourier transform of the function
$f(x),$ from ~(\ref{Q}), we obtain the integral form of the
remainder term for an arbitrary function $f(x)$ in the space
$L_{2}(0,b),$
$$
Q_{R}(f;x):=S_{R}(f(x);x)-f(x)=\frac{1}{2\pi
i}\oint_{|\lambda|=R}\frac{exp(i\lambda
x)}{\Delta(\lambda)}d\lambda(\frac{1}{2\pi}\int_{-\infty}^{\infty}\frac{\hat{f}(-\mu)\Delta(\mu)}{\mu-\lambda}d\mu).
$$

A similar formula was proved in another way by A.M. Sedletskiy in
the work [6]. Let us write the set of conditions
$$
\left\{%
\begin{array}{ll}
    \max supp\sigma(x)=b, &  \\
    \min supp(\sigma(x)+i)=0,&  \\
\end{array}
\right. \eqno (A)
$$
$$\sup|Im\lambda_{n}|=M<\infty, \eqno (B)$$
$$
\inf_{n\neq k}|\lambda_{n}-\lambda_{k}|>0, \eqno (C)
$$
$$
\sup m_{n}=m<\infty. \eqno (D)
$$
Denote $\omega(\lambda)=|\Delta(\lambda)|^{2}$. Let us write the
Muckenhoupt condition
$$
\sup_{l}(\frac{1}{|I|}\int_{I}\omega(\lambda)d\lambda)(\frac{1}{|I|}\int_{I}\omega^{-1}(\lambda)d\lambda)<\infty,
\eqno (E)
$$
where $I$ an arbitrary interval of the real axis. If the operator
$$
S^{+}:\sum_{\lambda_{n}\in\Lambda}P_{n}(x)exp(i\lambda_{n}x)\mapsto\sum_{Re\lambda_{n}>0}P_{n}(x)exp(i\lambda_{n}x)
$$
is bounded in the space $L_{2}(0,b)$, then we will say that the
basis $\{u_{k,n},k=0,...,m_{n}-1,\lambda_{n}\in \sigma(L)\}$ of
the space $L_{2}(0,b)$ is Riesz basis.

The following theorems are true.

{\bf Theorem 5.1} Let the conditions (A), (B), (C), (D) be valid
and assume that the function
$\omega(\lambda)=|\Delta(\lambda)|^{2}$ satisfies the Muckenhoupt
condition (E). Then the system
$\{u_{k,n},k=0,...,m_{n}-1,\lambda_{n}\in \sigma(L)\}$ is a Riesz
basis in the space $L_{2}(0,b)$.

{\bf Theorem 5.2}  Assume the conditions of Theorem 5.1, except,
perhaps, the condition (D). Then for every $f\in L_{2}(0,b)$, we
have
$$
\|x^{\frac{1}{2}}(b-x)^{\frac{1}{2}}Q_{r}(f,x)\|_{C[0,b]}\rightarrow0
$$
continuously for $r\rightarrow\infty$.

{\bf Theorem 5.3} Suppose the conditions of Theorem 5.1. Then for
every $f\in L_{2}(0,b)$ the coefficient sequence satisfies
$$
\{\{c_{k,n}\}_{k=0}^{\infty}\}_{n=0}^{\infty}\in l_{2},
$$
and $\|\{c_{k,n}\}\|_{l_{2}}\leq C(L)\|f\|_{L_{2}}.$

\vspace{4mm}

\begin{center}
{\bf REFERENCE}
\end{center}

[1] M.A. Naimark, Linear differential operators, GITTL, Moscow,
1954; English transl., Part II, Ungar, New York, 1968.

[2] M. Ruzhansky, V. Turunen, Pseudo-Differential Operators and
Symmetries, Birkhauser, 2010.

[3] B.E. Kanguzhin and M.A. Sadybekov, Differential operators on a
segment. Distributions of the eigenvalues, Gylym, Almaty, 1996.
[in russian].

[4] E.C. Titchmarsh, The zeros of certain integral functions //
Proc. London Math. Soc. 1926, V. 25, No 4, pp. 283--302.

[5] M.L. Cartwright, The zeros of certain integral functions //
The Quarterly Journal of Math., Oxford series, 1930, V. 1, No 1,
pp. 38--59.

[6] A.M. Sedletskiy, On biortogonal expansions of functions to
exponent series on a real segment // Usp. Mat. Nauk, 1982, V. 37,
No 5(227), pp. 51--95 [in russian].

\end{document}